\documentclass[a4paper, 12pt]{amsart}
\usepackage{eufrak, bbm, amsmath}
\usepackage{amssymb}

\numberwithin{equation}{section}

\newtheorem{theorem}{Theorem}[section]
\newtheorem{lemma}[theorem]{Lemma}
\newtheorem{corollary}[theorem]{Corollary}
\newtheorem{proposition}[theorem]{Proposition}
\newtheorem{definition}[theorem]{Definition}
\newtheorem{rema}[theorem]{Remark}
\newenvironment{remark}{\begin{rema} \em}{\end{rema}}

\def\Z{\mathbb{Z}}

\def\N{\mathbb{N}}
\def\Q{\mathbb{Q}}
\def\C{\mathbb{C}}
\def\P{\mathbb{P}}
\def\div{\mathrm{div}}
\def\O{\mathcal{O}}

\begin{document}

\title{Numerically Gorenstein surface singularities are 
  homeomorphic to Gorenstein ones}

\author{Patrick Popescu-Pampu }

\date{Appeared in Duke Math. Journal {\bf 159} no. 3 (2011), 539-559.}

\setcounter{section}{0}

\maketitle
\pagestyle{myheadings} \markboth{{\normalsize 
Patrick Popescu-Pampu}} 
{{\normalsize The topology of Gorenstein surface singularities}} 

\begin{abstract}
  Consider a normal complex analytic surface singularity. It is called
  \emph{Gorenstein} if the canonical line bundle is holomorphically
  trivial in some punctured neighborhood of the singular point 
  and \emph{numerically Gorenstein} if this line bundle is
  topologically trivial. The second notion depends only on the
  topological type of the singularity. Laufer proved in 1977 that, given a
  numerically Gorenstein  topological  type of singularity, \emph{every}
  analytical realization of it is Gorenstein if and only if one has
  either a Kleinian or a minimally elliptic topological type. The question 
  to know if \emph{any} numerically Gorenstein topology was
  realizable by \emph{some} Gorenstein singularity was left open. 
  We prove that \emph{this is indeed the case}. Our method is to 
  \emph{plumb holomorphically} meromorphic 2-forms obtained by adequate
  pull-backs of the natural holomorphic symplectic  
  forms on the total spaces of  the canonical line bundles of  complex curves. 
  More generally, we show that \emph{any normal surface singularity is 
  homeomorphic to a  $\Q$-Gorenstein 
  singularity whose index is equal to the smallest common denominator 
  of the coefficients of the canonical cycle of the starting singularity.} 
\end{abstract}

\section{Introduction}

Let $(X,0)$ be a normal complex analytic surface singularity. 
It  is called \emph{Gorenstein} if the
canonical line bundle of some representative of it is
holomorphically trivial outside $0$, and 
\emph{numerically Gorenstein} if the same line bundle is topologically
trivial outside $0$. 
For example, hypersurfaces and complete intersections with isolated
singularities  are normal and Gorenstein. 

 One has fundamental systems of neighborhoods 
of $0 \in X$ which are homeomorphic to cones over an oriented  3-manifold, 
called the \emph{boundary} or the \emph{link} of $(X,0)$. 
By a theorem of 
Neumann \cite{N 81}, the knowledge of the boundary up to oriented 
homeomorphisms is equivalent to the knowledge of the dual graph of 
the exceptional divisor of any good resolution, 
weighted by the genera and the 
self-intersections of its irreducible components. We say that two
singularities are  
\emph{homeomorphic} if and only if they have boundaries which are 
oriented homeomorphic. By Neumann's theorem, this means that they 
have good resolutions with isomorphic weighted 
dual graphs. 

The \emph{canonical cycle} of a resolution of the surface singularity is 
by definition the unique cycle supported by the exceptional divisor $E$,  
which has the same intersections with the irreducible components of $E$ 
as a canonical divisor. Its coefficients are determined by the weighted dual
graph, therefore they are topological invariants of the singularity. 
The coefficient of each component of $E$ is called its \emph{discrepancy}. 
The fact that a singularity is \emph{numerically Gorenstein} is a topological 
property: it is equivalent to the fact that, for  some 
resolution, all discrepancies are integral.  

Instead, the fact that a singularity is Gorenstein is \emph{not} in general
a topological  
property. In fact, Laufer \cite[Theorem 4.3]{L 77} proved that
\emph{all} analytical  
realizations of a given topology are Gorenstein if and only if one has 
either a Kleinian or a minimally elliptic topological type. 

It seems that the problem \emph{to decide if any numerically Gorenstein normal 
surface singularity is homeomorphic to a Gorenstein one} was left open. 
We learnt this problem from Jonathan Wahl (see the introduction 
of \cite{PS 09}).  We prove 
(Theorem \ref{numeqgor}):

\smallskip
  \emph{Any numerically Gorenstein normal surface singularity is homeomorphic 
  to a Gorenstein  one, having, moreover, an exceptional set of the minimal good
  resolution  analytically isomorphic to that associated to the initial singularity. }
\smallskip

One gets automatically a characterization of the weighted
graphs which  are realized by Gorenstein surface singularities: \emph{they are 
precisely the negative definite ones with  integral discrepancies}.
This fact was known before in the case of weighted homogeneous surface
singularities,   
as a consequence of Watanabe's criterion of \cite[(2.9)]{W 81} for such a
singularity to be  
Gorenstein (another proof of it is given by Wahl \cite[(4.4.2)]{W 88}).

Our proof starts from the fact that \emph{the holomorphic cotangent bundle of a
 complex manifold is  
naturally endowed with a holomorphic symplectic form}. We get in this way a
holomorphic    
symplectic form  on the canonical bundle $K_C= T^* C$ of any smooth, complex, 
compact curve  $C$.  We pull-back 
this form using meromorphic maps to the total spaces of convenient line
bundles. In order to take care of the components of $E$ with discrepancy $-1$, 
which are to be treated separately,   
we construct also special meromorphic 2-forms on the total spaces of
line  bundles on $\P^1$. 

The fundamental facts which make our constructions sufficient are the
following: 

$\bullet$  They lead to 
meromorphic 2-forms having \emph{local normal forms without moduli} in
holomorphic coordinates (namely, of the type $u^a v^b \: du \wedge dv$), which
allows to \emph{plumb} them \emph{holomorphically}.  

$\bullet$ If the initial numerically Gorenstein singularity is not Kleinian,
cusp or simple elliptic 
(cases when it is automatically Gorenstein), then all the components with
discrepancy $-1$ are rational and no two of them intersect. 

By plumbing tubular neighborhoods of the 0-sections \emph{using the
  identifications  given by the coordinate 
systems of the normal forms}, we get smooth,  
complex analytic surfaces endowed with a meromorphic 2-form which is 
\emph{symplectic} outside the compact divisor obtained as the image of
the 0-sections  
of the line bundles. Therefore, whenever the divisor is contractible,
one gets a normal Gorenstein singularity. 
 
 \medskip
 
 Teissier and Wahl  
 suggested we extend the theorem to all normal surface singularities,
 the role of Gorenstein singularities being played now by the
 $\Q$\emph{-Gorenstein} 
 ones (that is, those normal singularities such that 
 some power of the canonical bundle is holomorphically trivial in a pointed
 neighborhood of the singular point).  

We were able to prove this extension after Veys mentioned a result 
 of his (see Theorem \ref{core}), implying that  \emph{if the initial normal surface
 singularity is  not log canonical, then for the minimal good resolution all 
 components of $E$ with 
 discrepancy $-1$ are rational, they intersect at most two other components of $E$,
 and they are pairwise disjoint}.  

In order to state the extension, we need a little more vocabulary. The
\emph{index} of  
 a $\Q$\emph{-Gorenstein} singularity is the smallest exponent $r$ such 
 that $K_{X \setminus 0}^r$ is holomorphically trivial and the \emph{topological
   index} is defined in the same way for \emph{all} normal singularities, by
 asking $K_{X \setminus 0}^r$ to be \emph{topologically trivial}. By using the
 same kind of constructions as before, we prove
 (Theorem \ref{topqgor}):

\medskip
  \emph{Any normal surface singularity is homeomorphic 
  to a $\Q$-Gorenstein  one whose index is equal to the topological index 
  of the initial one, having, moreover, an exceptional set of the minimal good
  resolution  
  analytically isomorphic to that associated to the initial singularity. }
\medskip

 In this enlarged context, the analog of the natural holomorphic symplectic
 form on $K_C$ 
 is a natural meromorphic 2-form of weight $r$ on the total space of $K_C^r$,
 for all $r \in \N^*$. Namely, in natural local coordinates it is written
 $p^{1-r}\: (dp \wedge dq)^r$ (see Lemma \ref{gensym}). 

Wahl mentioned to us that he can prove the previous theorem in the weighted
homogeneous case by a generalization of Watanabe's arguments of \cite{W 81}
(unpublished).

 \medskip
 
 We describe now briefly the content of the paper. 
 In Section \ref{dual} we recall the fundamental notions about normal 
 surface singularities used in the rest of the paper.   In Section
 \ref{spemer} we  
 construct special meromorphic 2-forms on total spaces of line bundles on 
 smooth, compact, complex curves.  In Section \ref{mainthm} 
 we plumb such total spaces along weighted dual graphs,  
 in a way which leads to the proof of Theorem \ref{numeqgor}. Finally, 
 in Section \ref{exten}, we prove Theorem \ref{topqgor} by generalizing the
 steps of the proof of Theorem \ref{numeqgor}.

\section{Background} \label{dual}

In the sequel, the \emph{singularities} we consider will be \emph{germs of
  normal 
  complex analytic surfaces}. 
If $(X,0)$ denotes a singularity, tacitly $X$ will also denote a Milnor
representative  of it.

The divisor of a meromorphic $2$-form $\alpha$ on a smooth surface will
be denoted  $\mbox{div}(\alpha)$. Outside it, such a form is \emph{holomorphic 
symplectic}, that is, non-degenerate (as a consequence, by Darboux' theorem, it 
may be written $du \wedge dv$ in convenient local coordinates in the
neighborhood  
of any point, a fact we will not use in this paper). 

A resolution of a singularity is said to be 
\emph{good} if its exceptional divisor 
has normal crossings and its components are smooth. There exists a
\emph{minimal} good resolution,  
unique up to unique isomorphism. 
\medskip

\begin{definition}
  Let $(X,0)$ be a normal surface singularity. It is called \textbf{Gorenstein}
  if its  
  canonical bundle is holomorphically trivial in a punctured 
  neighborhood of $0$. It is called \textbf{numerically Gorenstein} 
  if one has instead only \emph{topological} triviality.
\end{definition}

In other words, a Gorenstein singularity is one which may be endowed 
with a 2-form which is holomorphic and non-vanishing in a punctured 
neighborhood of $0$, that is, which is 
holomorphic symplectic there.

Starting from a normal surface singularity, denote by 
$(\tilde{X}, E)\stackrel{\pi}{\rightarrow} (X,0)$ a good  resolution, 
where $E$ is the reduced exceptional divisor.  
By taking pull-backs by the morphism $\pi$, one identifies the space of
holomorphic  
2-forms on $X \setminus 0$ with the one of holomorphic 2-forms on 
$\tilde{X} \setminus E$. In particular, one may construct a holomorphic 
symplectic form on $X \setminus 0$ by constructing a holomorphic 
symplectic form on $\tilde{X} \setminus E$. It is the approach 
taken in this paper. 

Denote by $\Gamma$ the \emph{dual (intersection) graph} of $E$: its vertices
correspond bijectively to the irreducible components of $E$ and between two
distinct vertices $i$ and $j$
there are as many (unoriented) edges as the intersection number
$e_{ij}:=E_i \cdot E_j\geq 0$ of the corresponding components. In particular,
$\Gamma$ has no loops. Moreover, each vertex
$i$ of $\Gamma$ is weighted by
the self-intersection number $-e_i:=E_i^2 <0$
of the associated component $E_i$ inside $\tilde{X}$, 
and by the genus $g_i$ of $E_i$.

Denote by $V(\Gamma)$ the set of vertices of $\Gamma$. If $i\in V(\Gamma)$, 
its \emph{valency} is the number of edges  connecting it to other vertices.  
       
By classical results of 
Du Val \cite{V 44} and Mumford \cite{M 61}, the 
intersection form of the weighted dual graph 
is negative definite. Conversely, Grauert \cite{G 62} proved that if the
form associated to a reduced compact divisor $E$ on a
smooth surface is negative definite, then $E$  can be \emph{contracted}
to a normal singular point of an analytic surface, well-defined up to unique
analytic isomorphisms.

As the intersection form  is negative definite, there exists
a unique divisor \emph{with rational coefficients} $Z_{can}$ supported by
$E$ such that 
    $ Z_{can} \cdot E_i= K_{\tilde{X}} \cdot E_i$,  for all  $i \in V(\Gamma).$

\begin{definition}
The divisor $Z_{can}$ is called the {\bf canonical cycle} of $E$ (or of the
resolution $\pi$). 
If $Z_{can}= \sum_{i\in V(\Gamma)}k_i E_i$, then  
the coefficient $k_i$ is called the {\bf discrepancy} of $E_i$.  
\end{definition}

 The adjunction relations, applied to the curves $E_i$ inside the surface
 $\tilde{X}$, may be written as:
\begin{equation} \label{eqnum}
   (k_i + 1)e_i =  -2g_i + 2 + \sum_{\stackrel{j \in V(\Gamma)}{j \neq i}}
   k_j e_{ij}, \: \forall \: i \in V(\Gamma) .
 \end{equation}

Durfee \cite{D 78} proved the following characterization of numerically 
Gorenstein singularities (which explains the attribute ``\emph{numerically}''):

\begin{proposition} \label{durf}
   The normal surface singularity $(X,0)$ is numerically
    Gorenstein  if and only if for any (or some)
    good resolution, all discrepancies are integral.
\end{proposition}

One has the following generalization of the notion of normal Gorenstein
singularity: 

\begin{definition}
   The normal surface singularity $(X,0)$ is called {\bf $\Q$-Gorenstein} if
   there exists  
   a positive integer $r$ such that the pluricanonical bundle 
   $K^r_{X}$ is holomorphically trivial in a punctured neighborhood 
   of $0$. The smallest such exponent $r$ is then called the 
   {\bf index} of the singularity. 
\end{definition}

Equivalently, $(X,0)$ is $\Q$-Gorenstein if and only if the 
canonical class $K_X$ (consisting a priori of Weil divisors) 
is (holomorphically) $\Q$-Cartier and the index 
is the smallest integral factor making it Cartier. 
The normal Gorenstein singularities are precisely the
$\Q$-Gorenstein ones of index $1$. 

In dimension 2, being $\Q$-Gorenstein is equivalent to being a
cyclic quotient  
of a Gorenstein singularity (the \emph{canonical cover}), 
a fact which is no longer true in higher dimensions, where the  
canonical cover may not be Cohen-Macaulay. But in this paper we don't use this
view-point  
on $\Q$-Gorenstein surface singularities. 

If $(X,0)$ is $\Q$-Gorenstein and $K^r_{X}$ is holomorphically trivial 
on $X \setminus 0$, the pull-back by the resolution morphism $\pi$ of a
holomorphic trivializing section of $K^r_{X}$  
becomes a meromorphic section of $K^r_{\tilde{X}}$ whose divisor is exactly 
$rZ_{can}$. Therefore, in this case $r$ is a common denominator of the 
discrepancies. By the same method as the one used by Durfee to prove
Proposition \ref{durf}, one shows that the smallest common denominator of the
discrepancies is the smallest exponent $r\in \N^*$ such that $K^r_{X}$ is
topologically trivial on $X \setminus 0$. This motivates us to 
introduce the following:

 \begin{definition}
   The {\bf topological index} of the singularity $(X,0)$ is the smallest 
    positive integer $r$ such that 
   $K^r_{X\setminus 0}$ is topologically trivial.  
 \end{definition}

Recall also the following special classes of singularities:
\medskip

$\bullet$ the \emph{Kleinian singularities}, also known as
\emph{Du Val singularities, rational double points or simple surface
singularities} (see Durfee \cite{D 79}) are, up to isomorphism,
the surface singularities of the form $\mathbb C^2/G$, where $G$
is a finite subgroup of $SU(2)$. For these singularities, the
minimal resolution is good, it has $Z_K = 0$ and its dual graph $\Gamma$ is one 
of the trees $A_n, D_n, E_6, E_7, E_8$ with $g_i=0, e_i=2, \: \forall \: i\in
V(\Gamma)$.

\medskip

$\bullet$ the \emph{cusp singularities} are those normal surface
singularities which have a good resolution whose dual graph 
  $\Gamma$ is homeomorphic to a circle and $g_i=0, \:
 \forall \: i\in V(\Gamma)$   (see Looijenga \cite[page 16]{L 84}).
 
 \medskip

$\bullet$ the \emph{simple elliptic singularities}  are those normal surface
singularities whose minimal resolution has a smooth exceptional divisor 
of genus $1$. As well as the cusp singularities, they are \emph{minimally 
elliptic} (see Laufer \cite{L 77}). As we do not need this last notion
in the paper,  we do not give more details about it.

\medskip

$\bullet$ the \emph{log canonical singularities} are those normal 
surface singularities whose canonical cycle satisfies the inequalities 
$k_i \geq -1, \: \forall \:  i \in V(\Gamma)$.  By a result of Sakai
\cite{S 87},  \emph{they are necessarily $\Q$-Gorenstein, with an index equal
  to their  
topological index}. One may find their classification in Sakai \cite{S 87},  
Kawamata \cite{K 88},  Wahl \cite{W 90}, Alexeev \cite{A 92} or Matsuki
\cite{M 02}. By this classification,   
one sees that the numerically Gorenstein log canonical singularities (therefore 
the Gorenstein ones) are exactly the Kleinian, the cusp and the simple 
elliptic singularities. This class may also be defined in higher dimensions. 
Then, one usually asks the singularity to be $\Q$-Gorenstein, in order to 
be able to define the discrepancies (an exception is the recent paper 
\cite{BFF 10}, in which
discrepancies are defined also for singularities which are not
$\Q$-Gorenstein). Sakai's result quoted before shows that  
this supplementary condition is not needed in dimension 2 (see also 
\cite[Remark 4-6-3 (ii)]{M 02}). 
\medskip

In \cite[2.10]{V 04}, Veys proved the following theorem about the behaviour 
of the coefficients of the canonical cycle on the dual graph of the
\emph{minimal} good resolution:

\begin{theorem} \label{core}
    Let $(X,0)$ a normal surface singularity \emph{which is not log canonical}. 
    Denote by $V^{-}(\Gamma)$ the set of vertices $i$ of the minimal good 
    resolution such that $k_i < -1$, which is non-empty by hypothesis. 
    Then:
    
    \begin{enumerate}

   \item  \label{conn} 
      The full subgraph $\Gamma^{-}$ of $\Gamma$ with vertex set
      $V^{-}(\Gamma)$  is connected.

   \item \label{growing}  
      The full subgraph of $\Gamma$ with vertex-set $V(\Gamma) \setminus
      V^{-}(\Gamma)$  
      is a disjoint union of segments whose vertices represent rational
      components  
      of the exceptional divisor. Moreover, exactly one end of each segment 
      is attached (by one edge) to $\Gamma^{-}$, and the discrepancies are
      strictly increasing  
      along the segment, starting from that end. 
      
  \item \label{negat}
      All the discrepancies are (strictly) negative. 
     
     \end{enumerate}
\end{theorem}

In  \cite[2.11]{V 04}, Veys states that the previous theorem may be 
obtained as a corollary of results of Wahl \cite{W 90}. Instead, he uses
mainly results of Alexeev \cite{A 92}. 
Before, Laufer \cite[Prop. 2.1]{L 87} had proved the analog of 
(\ref{core}.\ref{negat}) for the {\it minimal} (not necessarily good) 
resolution of any normal surface singularity which is not Kleinian. 

\begin{remark}
 Veys' terminology is slightly different from ours. For
instance, he uses the name {\it log resolution} instead of  {\it good 
  resolution} and he works with the {\it log discrepancies} $a_i:= k_i +1$ 
instead of the discrepancies $k_i$. We prefer to work with the discrepancies
because 
in our computations they appear explicitely as exponents of coefficients of
meromorphic 2-forms (see the last two sections). 
\end{remark}

As a corollary of (\ref{core}.\ref{growing}), we have
the  following special case 
concerning the numerically Gorenstein singularities which are not log
canonical:

\begin{corollary} \label{sperel}
  Start from a numerically Gorenstein singularity which is not smooth, 
  Kleinian, a cusp singularity or a simple elliptic one. 
 Then the dual graph $\Gamma$ of the minimal good 
  resolution satisfies $k_i \leq -1$ for all $i \in V(\Gamma)$. Moreover,  
     if  $k_i = -1$, then $E_i$ is 
     a smooth rational curve, $i$ is
     of valency $1$ and the unique vertex $j$ to which it is connected
     satisfies $k_j = -2$. 
\end{corollary}

\noindent{\bf Proof.}  
   The only statement which is not a special case of the previous 
   theorem is the last one. But it is an immediate consequence 
   of the adjunction formulae (\ref{eqnum}).
   
   In fact the corollary could be proved by showing first that Laufer's 
   theorem \cite[Prop. 2.1]{L 87} stated before, namely, the fact that 
   (\ref{core}.\ref{negat}) holds for the minimal resolution of a
   normal surface 
   singularity which is not Kleinian, is also true for the minimal good
   resolution  
   of a numerically Gorenstein singularity. This is what we did before knowing 
   Veys' results. We analyzed  
   what happens at each step of the sequence of blowing-ups which 
   passes from the minimal resolution to the minimal good one and we 
   finished by using the adjunction relations. 
\hfill $\Box$

\medskip

One has a generalization for all normal surface singularities which are not
log canonical: 

\begin{corollary} \label{corgen}
   Start from a normal surface singularity which is not
   log canonical.  On the dual graph of the minimal good 
  resolution, 
     if $i$ is a vertex of $\Gamma$ such that $k_i = -1$, then it represents 
     a smooth rational curve and  its valency is $1$ or $2$. Moreover: 

     $\bullet$ \emph{If it is of valency} $1$, then the unique vertex $j$ to
     which it  is connected satisfies $k_j = -2$.

     $\bullet$ \emph{If it is of valency} $2$, then the two vertices $j_1,
     j_2$ to 
     which it is connected satisfy $k_{j_1} + k_{j_2} = -2$, with 
     $k_{j_1} \neq -1, \: k_{j_2} \neq -1$. 
\end{corollary}

\noindent{\bf Proof.} This is again a direct consequence of Theorem 
\ref{core} and of the adjunction formulae. 
\hfill $\Box$

\begin{remark} \label{adj}
    In particular, the  
    exceptional divisor of the minimal good resolution of non log canonical
    normal surface singularities  does not contain adjacent components $E_i$
    and $E_j$    such that $k_i= k_j=  - 1$ (see also Remark \ref{adjbis}). 
    As Veys told us, this conjectured statement was his motivation for proving
    Theorem \ref{core}.  
\end{remark}

\medskip
\section{Construction of special meromorphic 2-forms 
  on line bundles} \label{spemer}

Our idea for proving the theorem stated in the title of the paper, is to
construct first \emph{meromorphic} 2-forms on total spaces of 
line bundles which are \emph{holomorphic symplectic outside the union 
of the 0-section and of some fibers}, and then
to \emph{plumb} them in such a way that the symplectic forms glue.  This
section is dedicated to the first step and 
the next one to the second step. 

\medskip

When we say that $L$ \emph{is a line bundle}, $L$ will denote the total space
(and not the associated sheaf of regular sections). In particular, by \emph{a
  form 
defined on $L$} we will mean a form defined on the total space of $L$.

\begin{proposition} \label{consform}
  Let $C$ be a compact, connected and smooth complex curve of genus $g
  \geq 0$. 
  Fix an integral  divisor $D$ on $C$ and an integer $m \in \Z^*$ which
  divides $2g 
  -2 - \deg(D)$. Then there exists a line bundle $L
  \stackrel{\psi}{\longrightarrow} C$ of degree $m^{-1}(2g-2 -
  \deg(D))$, endowed with a meromorphic $2$-form $\alpha$ such that:
     
  \begin{enumerate}
     \item \label{divform}
       $\div(\alpha)= (m-1)C + \psi^*(D)$, where $C$
       is seen here as the $0$-section of $L$. 

     \item \label{prec} For each point $A\in C$, if $a \in \Z$
       denotes the multiplicity of $D$ at that point, one has a system
       of local coordinates $(q, f)$ centered at $A$, seen as a point of the
       $0$-section  
       of  $L$, such that $q$ is constant on the fibers of 
       $\pi$, that $f=0$ defines $C$ and moreover: 
         $$ \alpha= f^{m -1} q^{a} \: df \wedge  dq.$$
  \end{enumerate}

\end{proposition}

\noindent{\bf Proof.} Denote by $L_D$ the line bundle over $C$ whose
sheaf of sections  
is $\O(D)$ and by $s_D$ its canonical meromorphic section defining $D$. 
Tensoring by $s_D$ one gets a meromorphic map between total spaces of line 
bundles:
\begin{equation} \label{tens}
   H_D :=K_C \otimes L_{-D} \stackrel{\tau_D}{\cdots \longrightarrow} K_C, 
\end{equation}
where $K_C\stackrel{\gamma}{\longrightarrow} C$ denotes the canonical line
bundle of $C$, equal to its holomorphic  
cotangent bundle. 
Denote by $\alpha_C$ the canonical holomorphic symplectic form on the surface 
$K_C$. If $q$ denotes a local coordinate on the open set $U \subset C$ and $p$
is the associated  
dual coordinate on $\gamma^{-1}(U)\subset K_C$ (that is,
the  general   
$1$-form on $U$ is written $pdq$), one has the following expression for 
$\alpha_C$ on $\gamma^{-1}(U)$ :
  \begin{equation} \label{cansym}
      \alpha_C= dp \wedge dq.
  \end{equation}
Denote then by $\alpha_D$ the meromorphic 2-form on $H_D$
obtained as the pull-back of $\alpha_C$ by $\tau_D$:
  \begin{equation} \label{preim1}
     \alpha_D:= \tau_D^*(\alpha_C).
  \end{equation} 
 
 Consider now any $m$-th root $H_{D}^{1/m}$ of the line bundle $H_D$ 
 and the $m$-th tensor power meromorphic map:
 \begin{equation} \label{mtens}
   H_{D}^{1/m} \stackrel{\mu_{D,m}}{\cdots \longrightarrow} H_D. 
\end{equation}
Such line bundles $H_{D}^{1/m}$ are the solutions $H$ of the equation $m\cdot H =
H_D$ in   
the Picard group $(\mathrm{Pic}(C), +)$ written additively. The set of 
solutions is non-empty by the hypothesis that $m \neq 0$ divides 
$\deg(H_D)= 2g-2 - \deg(D)$ and it is a principal 
space (a torsor) under the subgroup of $|m|$-torsion elements of the Jacobian
$\mathrm{Pic}^0(C)$ of $C$. This Jacobian is a real $(2g)$-dimensional torus, 
therefore there are $|m|^{2g}$ such line bundles. 
We may take any one of them.

Denote by $\alpha_{D,m}$ the meromorphic 2-form on $H_{D}^{1/m}$
obtained as the pull-back of $\alpha_D$ by $\mu_{D,m}$:
  \begin{equation} \label{preim2}
     \alpha_{D,m}:= \mu_{D,m}^*(\alpha_D).
  \end{equation}

  We claim that 
  \emph{the pair $(L, \alpha): = (H_{D}^{1/m}, \alpha_{D,m})$
   satisfies both   
  properties \ref{divform}. and \ref{prec}}. 
\medskip  

  The rest of the proof
  is devoted to an  explanation of this claim. 
  
  The map $\tau_D\circ \mu_{D,m}$ 
  (see (\ref{tens}) and (\ref{mtens})) 
  is regular and without critical 
    points on the complement $H_{D}^{1/m}\setminus (C \cup
    \psi_{D,m}^{-1}(|D|))$, where  
    $H_{D}^{1/m} \stackrel{\psi_{D,m}}{\longrightarrow} C$ denotes the
    bundle projection and $|D|$ denotes the support of $D$.  
    It is therefore enough to check point {\it \ref{prec}}. in order to deduce
    also point {\it \ref{divform}}.

  We will proceed in two steps, by computing first the form $\alpha_D$ 
  in adapted local coordinates and then the form $\alpha_{D,m}$. 
  
  \medskip
  $\bullet$   Consider any local coordinate $q$ on $C$ centered at
  $A$. Choose then a local  
  trivializing section $\sigma_A$ of $L_D$ in the neighborhood of $A$
  such that  
  $s_D= q^a \sigma_A$. Denote by $(\sigma_A)^*$ the dual trivializing
  section of  
  $L_{-D} = L_D^*$.  Then $dq \otimes (\sigma_A)^*$ is a local
  trivializing section  
  of $H_D$, with respect to which a general
  section is written  
  $h \: dq \otimes (\sigma_A)^*$. Then $(q,h)$ is a system of local
  coordinates on the  
  total space $H_D$ in the neighborhood of $A$. 
  
  We have $\tau_D(h \: dq \otimes (\sigma_A)^*)= h \: dq \:
  s_D((\sigma_A)^*)= h \: q^a dq$,  
  which shows that with respect to the local coordinate systems
  $(q,h)$ on $H_D$ and  
  $(q,p)$ on $K_C$, the map $\tau_D$ is given by: $(q,h) \rightarrow
  (q, hq^a)$.  
  Combined with formula (\ref{cansym}), this implies that the meromorphic 
  2-form $\alpha_D$  
  defined by (\ref{preim1}) has the following 
  expression in the coordinate system $(q,h)$ on $H_D$:
    \begin{equation} \label{exprad1}
       \alpha_D= q^a \: dh \wedge dq.
    \end{equation}
    
    $\bullet$ We do now the same kind of reasoning for the map 
       $\mu_{D,m}$ (see (\ref{mtens})). Choose a local trivializing
       section $\eta_A$ of  
       $H_{D}^{1/m}$ in the neighborhood of $A$, 
       such that $\mu_{D,m}(\eta_A)= \eta_A^{\otimes m} = dq \otimes
       (\sigma_A)^*$.  
       If a general section of $H_{D}^{1/m}$ is written $ l \:
       \eta_A$, then  $(q,l)$ is 
       a local coordinate system on the total space $H_{D}^{1/m}$ in
       the neighborhood  
       of the point $A$. With respect to the local coordinate systems
       $(q,l)$ on  
       $H_D^{1/m}$ and $(q,h)$ on $H_D$, the map $\mu_{D,m}$ is given by 
       $(q,l) \rightarrow (q, l^m)$. 
        Combined with formula (\ref{exprad1}), this implies that the
        meromorphic 2-form $\alpha_{D,m}$  
  defined by (\ref{preim2}) has the following 
  expression in the coordinate system $(q,l)$ on $H_D^{1/m}$:
    \begin{equation} \label{exprad2}
       \alpha_{D,m}= m q^a l^{m-1} \:  d l  \wedge dq.
    \end{equation}

By making the change of variables $f := m^{1/m} \: l$,
where $m^{1/m}$ is an arbitrary  
$m$-th root of $m$ (which exists, by the hypothesis $m \neq 0$), we get from
equation (\ref{exprad2}): 
  $$  \alpha_{D,m}= f^{m -1}q^{a}  \: df \wedge dq.$$
  
  The proposition is proved.
\hfill $\Box$

\begin{remark} The statement (\ref{consform}.{\it \ref{prec}.}) is
  not a consequence  of the statement (\ref{consform}.{\it \ref{divform}}.). 
  Indeed, there are moduli for germs of  
meromorphic 2-forms whose divisor has normal crossings with two
irreducible components, implying that not all germs of meromorphic 2-forms whose
divisor has normal crossings may be written in the form $u^a v^b \: du \wedge
dv$ in adapted local coordinates.  One gets the simple normal forms stated in
(\ref{consform}.{\it \ref{prec}.}) as a consequence of our construction  by
successive pull-backs  
of the standard symplectic form on the canonical bundle.
\end{remark}

The previous proposition does not allow to get 
meromorphic $2$-forms having a pole of order $1$ along the $0$-section. 
In the next one, we construct special meromorphic 2-forms on the
total spaces  
of arbitrary line bundles on a smooth \emph{rational} curve, having 
such pole of order $1$ along the $0$-section. We will see  
that this is enough to prove the characterization of the topological types of 
Gorenstein singularities. 

\begin{proposition} \label{ratcomp}
   Let $C$ be a compact, connected and smooth rational curve. Let
   $L\stackrel{\psi}{\longrightarrow} C$ 
   be a line bundle on it  and $A$ be a point on
   $C$. Then there exists a meromorphic 2-form $\alpha$ on $L$ such that: 
  \begin{enumerate}
     \item $\div(\alpha)= -C - 2 \: \psi^*(A)$; 

     \item   one has a system
       of local coordinates $(q, f)$ centered at $A$, seen as a point of the 
       $0$-section of $L$, such that $q$ is constant on the fibers of
       $\psi$, that $f=0$ defines $C$ and moreover: 
         $$ \alpha= f^{-1} q^{-2} \:  df \wedge dq.$$
  \end{enumerate}
\end{proposition}

\noindent{\bf Proof.} Choose a covering of $C$ by two affine charts
isomorphic to $\C$,  
one with coordinate $q_1$ centered at $A$ and the other with coordinate
$q_2$, such that:  
\begin{equation} \label{inv}
     q_2=q_1^{-1}. 
\end{equation}
Choose then a trivialization of $L$ over the $q_1$-chart, with respect to which 
the fiber coordinate is denoted $f_1$ and another trivialization over
the $q_2$-chart,  
with respect to which the fiber-coordinate is denoted $f_2$, such that
the change of trivialization over the overlap is given by the equation:
  \begin{equation} \label{fiber}
     f_2= q_1^n f_1.
  \end{equation} 
Here $n= \deg(L)$. 
Such trivializations exist because on $\P^1$ any line bundle of degree $n$ is 
isomorphic to $\O(n)$. Therefore:
   \begin{equation} \label{form}
     df_2 \wedge dq_2= - q_1^{n-2} df_1 \wedge dq_1.
  \end{equation} 
Using also the equations (\ref{inv}) and (\ref{fiber}), this implies that the
meromorphic 2-form which on the chart with  coordinates 
$(q_2,f_2)$ of $L$  is given by 
  $\alpha:= - f_2^{-1} \: df_2 \wedge dq_2$, 
becomes 
 $\alpha= f_1^{-1} q_1^{-2}  \:  df_1 \wedge dq_1$ 
in the chart with coordinates $(q_1,f_1)$.
  Looking at the expressions of $\alpha$ in both charts, we see that the
  proposition is proved by taking $(q,f):=(q_1,f_1)$. 
\hfill $\Box$

\medskip
\section{The topology of Gorenstein normal surface
  singularities} \label{mainthm} 

We are ready to establish the theorem announced in the title of the
paper. The idea is to \emph{plumb} using biholomorphic maps the line
bundles constructed using Proposition \ref{consform}  
and Proposition \ref{ratcomp}, in such a way that the associated
meromorphic 2-forms \emph{glue}. We recall that the
operation of \emph{plumbing} of tubular neighborhoods of $0$-sections of line
bundles over compact Riemann surfaces was introduced  by Mumford in
\cite{M 61}.  Since then, this operation 
was nevertheless mainly used in \emph{the differentiable category}, that is, 
without necessarily respecting holomorphic structures. Instead, in our 
approach it is essential to do \emph{holomorphic} plumbing. This procedure 
was also used by Laufer (see e.g. \cite[page 83]{L 73}).

\begin{theorem} \label{numeqgor}
  Consider a numerically Gorenstein surface singularity. Then it is
  homeomorphic to a normal Gorenstein surface singularity. Moreover, 
  the second singularity may be chosen such that the exceptional curve of 
  its minimal good resolution is analytically isomorphic to that of the initial 
  singularity. 
\end{theorem}

\noindent{\bf Proof.} 
Let $(X,0)$ be a \emph{numerically Gorenstein} singularity 
and $(\tilde{X}, E) \stackrel{\pi}{\longrightarrow} (X,0)$ its minimal 
good resolution. 
If the singularity is smooth, Kleinian, a cusp singularity or a simple
elliptic one,  then it is automatically Gorenstein (see Laufer \cite{L 77}). 
Therefore, we may suppose that \emph{we are not in 
one of these cases}. Consequently, we may apply Corollary \ref{sperel}. 

For each vertex $i$ of $\Gamma$, denote by $V_i\subset V(\Gamma)$ the set of 
indices of those components of $E$ which have non-empty intersection with
$E_i$.  
For every $j \in V_i$, denote $E_{i,j} := E_i \cap E_j $. This set may 
consist of several  points: $E_{i,j}= \{A_{i,j,l} \: | \: l \in I(i,j)\}$. Of
course, $E_{i,j} = E_{j,i}$, which allows us to choose the index sets such that 
$A_{i,j,l} = A_{j,i,l}$ for all $l \in I(i,j) = I(j,i)$. We will look at  
$E_{i,j}$ also as a reduced divisor on the abstract curve $E_i$. 

By Corollary \ref{sperel}, we have $k_i \leq - 1$ for
all $i \in I$.  We consider two cases:

$\bullet$  \emph{Suppose  that $k_i < -1$.} Then apply Proposition 
\ref{consform} to $C:= E_i$,  
  $g:= g_i$, $m:= k_i +1$, $D:= \sum_{j \in V_i} k_j E_{i,j} $.   The
  hypothesis that $m$ divides $2g - 2 - \deg(D)$ is a consequence of formula
  (\ref{eqnum}).  

$\bullet$  \emph{Suppose  that $k_i = -1$}. By Corollary \ref{sperel},  
we know that $E_i$ is then a smooth \emph{rational} curve and that $V_i$ has
only  
one element $j$, which satisfies moreover $k_j= -2$. 
Then apply Proposition \ref{ratcomp} to $C:= E_i$, $L:=\mathcal{O}(-e_i)$,
$\{A\}= E_i \cap E_j$. 
\medskip

By combining both cases, we see that for \emph{all} $i \in V(\Gamma)$, there
exists a line bundle  
$L_i\stackrel{\psi_i}{\longrightarrow} E_i$ \emph{of degree} $-e_i$ on $E_i$
and  a meromorphic 2-form $\alpha_i$ on it, such that:

     $\bullet$ $\div(\alpha_i)=  k_i \: E_i  +  \sum_{j \in V_i} k_j \:
       \psi_i^*(E_{i,j})$;  

     $\bullet$ for each point $A_{i,j,l}\in E_{i,j}$,  one has a system
       of local holomorphic coordinates $(q_{i,j,l}, f_{i,j,l})$ centered at
       $A_{i,j,l}$  on the total 
       space $L_i$, such that $q_{i,j,l}$ is constant on the fibers of
       $\psi_i$, that $f_{i,j,l}=0$ defines the $0$-section $E_i \subset L_i$ 
       and moreover: 
        \begin{equation} \label{formnorm}
             \alpha_i= (q_{i,j,l})^{k_j} (f_{i,j,l})^{ k_i} \:  df_{i,j,l} \wedge
             dq_{i,j,l}.  
         \end{equation}

Consider then neighborhoods  of  
the $0$-sections $E_i \subset L_i$,  and
\emph{plumb} them by using the  
following \emph{holomorphic} identifications of neighborhoods of the points
$A_{i,j,l}$ and $A_{i,j,l}$ contained in the previous charts,  whenever $i$
and $j$ are adjacent   
vertices of $\Gamma$ and $l \in I(i,j)$:
  $$\left\{ \begin{array}{ll} 
           q_{i,j,l}= \lambda_{i,j} f_{j,i,l} \\
            q_{j,i,l} = \lambda_{j,i} f_{i,j,l}
       \end{array} \right. \mbox{ with } \lambda_{i,j},\lambda_{j,i} \in \C^*.$$
 The non-vanishing complex numbers $\lambda_{i,j}, \lambda_{j,i}$ are
 to be chosen  
 such that by this change of coordinates, the form $\alpha_i$ given by
 (\ref{formnorm})  is transformed into $\alpha_j$, which says 
 precisely that 
\emph{they glue on the surface obtained by the 
 plumbing}. 

These compatibility conditions translate into the 
 following system of  equations, which is non-trivial because, when we plumb,
 we \emph{permute} the order of the coordinates: 
   \begin{equation} \label{coef}
        \lambda_{i,j}^{k_j +1} =- \lambda_{j,i}^{k_i  + 1} \:
         \mbox{ for all pairs } (i,j) 
         \mbox{ of adjacent vertices of the graph } \Gamma. 
    \end{equation}

As before, our hypothesis allows us to apply Corollary  \ref{sperel},
which implies  
that there is no pair of adjacent vertices such that $k_i =k_j=  -1$
(otherwise, equation  
(\ref{coef}) would be contradictory). In fact this is the only
obstruction to the existence of solutions   
of the system (\ref{coef}): orient each edge of $\Gamma$ such that, whenever 
an edge goes from a vertex $i$ to a vertex $j$, we have 
$k_j  < -1$. Define then  
$\lambda_{j,i}:=1$ and $\lambda_{i,j}$ to be a $(k_j + 1)$-th root of
$-1$. We get obviously  
a solution of (\ref{coef}).

After such a choice, the surface $\tilde{Y}$ obtained by the plumbing
is smooth complex analytic and  the 
images of the $0$-sections $E_i \subset L_i$ form a divisor canonically 
identified to $E$ and with the given weighted dual graph. 
By Grauert's theorem (see \cite{G 62}), 
this divisor $E$ may be contracted in $\tilde{Y}$, yielding a
normal singularity $(Y,0)$.  
By our construction, 
\emph{the meromorphic 2-forms $\alpha_i$ glue into a meromorphic 2-form
$\tilde{\alpha}$  
on $\tilde{Y}$, which is symplectic outside $E$}. Therefore, it descends into 
a holomorphic symplectic form on $Y  \setminus 0$,  which means 
that $(Y,0)$ is Gorenstein. 

The second statement of the theorem is a consequence of the fact that, 
in our construction, one starts only from an abstract compact curve
  $E= \cup_{i\in V(\Gamma)}E_i$ with simple  
  normal crossings, decorated with weights $e_i$ such that
  the associated  intersection form is negative definite and such that the
  coefficients $k_i$  determined by the adjunction system (\ref{eqnum}) are
  integers.  
  As a result of our  
  construction, one gets a Gorenstein singularity having a resolution whose 
  exceptional divisor is $E$,  and such that $E_i^2 = -e_i$ for each
  irreducible  component $E_i$.

The theorem is proved.
\hfill $\Box$

\begin{remark} \label{adjbis}
    We could do the plumbing construction using only the meromorphic 
    2-forms described in  
   Propositions \ref{consform} and \ref{ratcomp},  because we knew that we
   could work with a good resolution such that for any 
   two intersecting components $E_i, E_j$ of the exceptional divisor, 
   one does not have simultaneously $k_i = k_j =  -1$ 
   (see also Remark \ref{adj}). Indeed, in Proposition \ref{consform} 
   one has $m-1  \neq -1$ and in Proposition \ref{ratcomp}  one has adjacent 
   components of the divisor of the meromorphic 2-form with orders $1$ and
   $2$ (compare to Corollary \ref{sperel}).  
\end{remark}

\begin{remark} \label{spesing}
      Our construction yields singularities whose minimal good resolution has
      the following special property: \emph{each irreducible component of the
        exceptional divisor has a neighborhood isomorphic to a neighborhood of
      the 0-section of its normal bundle}. 
\end{remark}

The following corollary settles a question left open by Neumann and Wahl \cite{NW 05}: 

\begin{corollary}
  Any normal surface singularity whose boundary is an integral homology 
  sphere is homeomorphic to a Gorenstein normal surface singularity.
\end{corollary}

\noindent{\bf Proof.} 
  The boundary is an integral homology sphere if and only if $\Gamma$ is a tree 
  with $g_i=0$ for all $i \in V(\Gamma)$ and the determinant of the
  intersection form  
  is $\pm 1$. Therefore, the adjunction system (\ref{eqnum}) shows that, in
  this case,  the discrepancies 
  are all integral. We conclude using Theorem \ref{numeqgor}.
\hfill $\Box$

\medskip
\section{Extension to $\Q$-Gorenstein singularities} \label{exten}

If $S$ is a smooth complex manifold of dimension $s >0$ 
and $r \in \Z^*$, we will say that the holomorphic (meromorphic)
sections of the {\it pluricanonical line bundle} $K_S^r$ {\it of exponent $r$}
are 
{\it holomorphic (meromorphic) $s$-forms of weight $r$}.  If $(q_1,...,q_s)$ 
is a system of local coordinates on $S$, such a form is 
written canonically as :
  $$     f(q_1,...,q_s) (dq_1 \wedge \cdots \wedge dq_s)^r, $$
  where $f$ is holomorphic (meromorphic). 
  
  In particular, consider $S=C$, a smooth holomorphic curve (a Riemann 
  surface). Denote by: 
      $K_C^r \stackrel{\gamma_r}{\longrightarrow} C$
  the projection morphism.  
  If $q$ is a local coordinate on the open 
  set $U \subset C$, then $(dq)^r$ is a 
  local section of  $K_C^r |_U$. One may take 
  as a second coordinate on $\gamma_r^{-1} (U)$ the analog of the dual
  coordinate $p$  
  of the case $r=1$ (see the beginning of the proof of Proposition
  \ref{consform}), that is,  the holomorphic function 
  $p_r$ such that the sections $\alpha$ of  $K_C^r |_U$ (that is, the 1-forms
  of weight $r$) are written as:
      $$\alpha = p_r  (dq)^r. $$

The next lemma shows that \emph{all} pluricanonical bundles over a curve
are endowed with 
natural \emph{meromorphic} 2-forms of the same weight as the exponent of the
bundle, generalizing the existence of the natural symplectic form on the
canonical bundle.

\begin{lemma} \label{gensym}
   Let $C$ be a  smooth holomorphic curve and $r \in \Z^*$. Then the
   holomorphic $2$-forms  
   of weight $r$  on the total space $K_C^r$ which are expressed in local
   coordinates as follows:
   \begin{equation} \label{canplur2}
        p_r^{1 - r} \: (dp_r \wedge d q)^r
   \end{equation}
   glue into a well-defined meromorphic $2$-form $\alpha_{C,r}$ of weight $r$
   on the total space $K_C^r$. 
\end{lemma}
 
 \noindent{\bf Proof.} 
        Let $\tilde{q}$ be a second local coordinate on $C$, defined on 
        an open set overlapping the domain of definition of $q$. Denote 
        by $\tilde{p}_r$ the associated complementary coordinate on the total 
        space $K_C^r$. Then $\tilde{p}_r (d \tilde{q})^r = p_r  (dq)^r$, 
        which implies:
           \begin{equation} \label{canplur3}
               \tilde{p}_r =  p_r \cdot (\frac{d \tilde{q} }{d q} )^{ -r}.
           \end{equation} 
         As a consequence, 
           $d \tilde{p}_r \wedge d \tilde{q} =  (\frac{d \tilde{q} }{d q} )^{
             -r +1} dp_r   \wedge dq$
         which implies, by taking $r$-th powers on both sides and by applying
         again  
         relation (\ref{canplur3}), the desired equality:
           $$  \tilde{p}_r^{1-r} ( d \tilde{p}_r \wedge d \tilde{q})^r = 
                  p_r^{1 - r} \: (dp_r \wedge d q)^r. $$
       \hfill $\Box$

 The following proposition generalizes Proposition \ref{consform}. 
 It specializes to it for $r=1$. 
 
 \begin{proposition}  \label{consform2}
      Let $C$ be a compact, connected and smooth complex curve of genus $g
  \geq 0$. 
  Fix an integral divisor $D$ on $C$, an integer $r \in \N^*$ 
  and an integer $m \in \Z^*$ which divides $r(2g -2) - \deg(D)$. 
  Then there exists a line bundle $L
  \stackrel{\psi}{\longrightarrow} C$ of degree $m^{-1}(r(2g-2) - \deg(D))$, 
  endowed with a meromorphic $2$-form $\alpha$ of weight $r$, such that:
  \begin{enumerate}
     \item \label{divform2}
       $\div(\alpha)= (m-r)C + \psi^*(D)$, where $C$
       is seen here as the $0$-section of $L$. 

     \item \label{prec2} For each point $A\in C$, if $a \in \Z$
       denotes the multiplicity of $D$ at that point, one has a system
       of local coordinates $(q, f)$ centered at $A$, seen as a point of the 
       $0$-section of  $L$,  
       such that $q$ is constant on the fibers of 
       $\psi$, that $f=0$ defines $C$ and moreover: 
         $$ \alpha= f^{m -r} q^a  \: (df \wedge  dq)^r.$$
  \end{enumerate}
     
 \end{proposition}

  \noindent{\bf Proof.} 
       We follow the steps of the proof of Proposition \ref{consform}. 
       Consider the natural meromorphic map got by tensoring with the
       meromorphic  section $s_D$ of $L_D$:
          $$ H_{D,r} := K_C^r \otimes L_{-D} 
              \stackrel{\tau_{D,r}}{\cdots \longrightarrow} K_C^r.$$
        Then $(dq)^r \otimes \sigma_A^*$ is a trivializing section of 
        $K_C^r \otimes L_{-D}$. If $h_r (dq)^r \otimes \sigma_A^*$ denotes 
        the general section of $ H_{D,r}$, then $(q, h_r)$ is a local coordinate
        system on $ H_{D,r}$ in the neighborhood of $A$. With respect to the
        coordinate systems $(q, h_r)$ and  
        $(q, p_r)$, the map $\tau_{D,r}$ is given by 
          $ (q, h_r) \rightarrow (q, h_r q^a).$ 
        Recalling that $\alpha_{C,r}$ was defined by Lemma \ref{gensym}, this
        implies that:
          $$ 
                    (\tau_{D,r})^* \alpha_{C,r} 
                     = (\tau_{D,r})^* (p_r^{1-r} \: (d p_r \wedge d q)^r) =
                          (h_r q^a)^{1-r}\: (d(h_r q^a) \wedge dq)^r  = 
                          h_r^{1-r} q^a  \: (dh_r \wedge dq)^r.
                  $$
                 
           Consider now an arbitrary $m$-th root $H_{D,r}^{1/m}$ of the line 
           bundle  $H_{D,r}$, which exists by the hypothesis that $m\neq 0$
           divides $\deg(H_{D,r})= r(2g-2) - \deg(D)$. In adapted local
           coordinates, chosen as in the 
           proof  
           of Proposition \ref{consform}, the $m$-th tensor power map 
              $$ H_{D,r}^{1/m}  \stackrel{\mu_{D,r,m}}{\cdots
                \longrightarrow} H_{D,r}$$ 
         is given by 
             $ (q, l) \rightarrow (q, l^m). $
          Therefore:
                  $$ 
                    \alpha_{D,r,m} := (\mu_{D,r,m})^* (\tau_{D,r})^*
                    \alpha_{C,r} 
                          = l^{m(1-r)}q^a \: (m l^{m-1} \: d l 
                           \wedge dq)^r  = 
                        m^r l^{m-r}q^a  \: (dl \wedge dq)^r.
                   $$
           By making the change of variables $f := m^{r/m}l$, where $m^{1/m}$
           is an arbitrary $m$-th root of $m$, we get: 
             $$\alpha_{D,r,m} = f^{m-r} q^a\:  (df \wedge dq)^r.$$ 
  This shows that the pair $(L, \alpha) := (H_{D,r}^{1/m}, \alpha_{D,r,m})$
  satisfies both 
  properties {\it \ref{divform2}.} and {\it \ref{prec2}}. 
       \hfill $\Box$  
        \medskip

 In the same way as Proposition \ref{consform} did not allow to get
 meromorphic 2-forms with poles of order 1 along the 0-section, the previous
 one does not allow to get meromorphic 2-forms of weight $r$ with poles of
 order $r$ along 
 the 0-section. We get such poles using  
 the following proposition, which  generalizes Proposition \ref{ratcomp}. It 
 specializes to it for $r=1, m_1 = -2$ and  $m_2 = 0$. 
 
 \begin{proposition} \label{ratcomp2}
      Let $C$ be a compact, connected and smooth rational curve. Let 
   $L \stackrel{\psi}{\rightarrow} C$
   be a line bundle on it and $A_1, A_2$ be two distinct
   points on 
   $C$. Fix an integer $r \in \N^*$ and two other integers  
   $m_1, m_2 \in \Z \setminus \{ -r \}$ 
   such that $m_1 + m_2 = -2r$. 
   Then there exists a meromorphic 2-form $\alpha$ of weight $r$ on $L$ 
   such that: 
  \begin{enumerate}
     \item $\div(\alpha)= - r \: C  + m_1 \: \psi^*(A_1)  + m_2 \:
       \psi^*(A_2)$.  

     \item  \label{spefor} For each $l \in \{1,2 \}$, one has a system
       of local coordinates $(u_l, v_l)$ centered at $A_l$, seen as a point of
       the  $0$-section of  $L$,  such that
       $u_l$ is constant on the fibers of 
       $\psi$, that $v_l =0$ defines $C$ and moreover: 
         $$ \alpha= v_l^{-r} u_l^{m_l}\:  (dv_l \wedge du_l)^r.$$
  \end{enumerate}
 \end{proposition}
 
  \noindent{\bf Proof.} 
     We follow the notations introduced in the proof of 
  Proposition \ref{ratcomp}. That is, we work inside the charts 
  $(q_1,f_1), (q_2, f_2)$, chosen moreover such that $A_l$ is defined by 
  $q_l = f_l =0$. Consider the following meromorphic 2-form of weight
  $r$    on $L$:
   $$ \alpha := f_2^{-r} q_2^{m_2}  \:  (df_2 \wedge dq_2)^r.$$
  The equations (\ref{inv}), (\ref{fiber}), (\ref{form}) imply 
  immediately that, written in the other chart: 
     $$\alpha = (-1)^r f_1^{-r} q_1^{m_1} \:  (df_1 \wedge dq_1)^r. $$
 We almost get the desired theorem, with the exception of the fact that 
 in this second chart, one has the extra factor $(-1)^r$. We have to 
care about it, in order to get the same shape for \emph{all} the 
local normal forms, adapted to the plumbing operations to be done in the
sequel (in the proof of Theorem \ref{topqgor}). 

$\bullet$ \emph{If $r$ is even}, there is nothing more to be done. 

$\bullet$ \emph{If $r$ is odd}, then we do the change of variable  
  $q_1= \nu \cdot u_1$, where $\nu$ is an arbitrary $(m_1 + r)$-th root of
  $(-1)^r$, which exists, by the hypothesis $m_1 \neq -r$. We get: $\alpha =
  f_1^{-r} u_1^{m_1} \:  (df_1 \wedge du_1)^r.$  
  By taking $u_2 := q_2$ and $v_l := f_l$ for $l \in \{1,2 \}$, we see 
 that $\alpha$ is expressed as in point {\it \ref{spefor}.} in both charts 
 with coordinates $(u_l, v_l)$. As these charts cover $L$, the proposition is
 proved. 
            \hfill $\Box$
\medskip

 Now we are ready to prove a generalization of Theorem \ref{numeqgor} 
 to all normal surface singularities. It specializes to it for $r=1$. 
 
 \begin{theorem} \label{topqgor}
     Consider a normal surface singularity, with topological index $r \geq 1$. 
     Then it is homeomorphic to a $\Q$-Gorenstein singularity with index $r$. 
     In particular, being $\Q$-Gorenstein imposes no topological restrictions
     on normal surface singularities. Moreover, 
  the second singularity may be chosen such that the exceptional curve of 
  its minimal good resolution is analytically isomorphic to that of the initial 
  singularity. 

 \end{theorem}
 
  \noindent{\bf Proof.} The proof follows the same steps as the one of Theorem
  \ref{numeqgor} and it uses the same notations. For each vertex $i \in
  \Gamma$, we consider a convenient line bundle endowed with a meromorphic
  2-form of 
  weight $r$, constructed either by Proposition \ref{consform2}, or by
  Proposition \ref{ratcomp2}. We consider three cases,
  according to the values of the discrepancy $k_i$ and of the valency
  $k_i$. When   $k_i   = -1$, 
  we know by Corollary \ref{corgen} that the valency $v_i$ is either $1$ or
  $2$ and that $E_i$ is a rational curve.

 $\bullet$ \emph{Suppose that} $k_i \neq -1$. Then apply Proposition
 \ref{consform2} to $C:= E_i$, $g:= g_i$, $m:= r(k_i + 1)$, $D:= \sum_{j \in
   V_i} r \cdot k_j \: E_{i,j}$. The hypothesis that $m$ divides $r(2g-2) -
 \deg(D)$ is a consequence of formula (\ref{eqnum}). 

 $\bullet$ \emph{Suppose that $k_i = -1$ and $v_i =1$}. Denote by $j$ the
 unique element of $V_i$. Then apply Proposition \ref{ratcomp2} to $C:= E_i$,
 $L := 
 \mathcal{O}(-e_i)$, $\{ A_1\}= E_i \cap E_j$, $m_1 = -2r$ and $m_2= 0$. 

 $\bullet$ \emph{Suppose that $k_i = -1$ and $v_i =2$}. Corollary \ref{corgen}
 tells us that, if $j_1 \neq j_2$ are the two elements of $V_i$, then $k_{j_1}
 + k_{j_2} = -2$ and $ k_{j_1} \neq -1$, $ k_{j_2} \neq -1$. Then apply
 Proposition  \ref{ratcomp2} to $C:= E_i$, $L := 
 \mathcal{O}(-e_i)$, $\{ A_l\}= E_i \cap E_{j_l}$ and  $m_l = r \cdot
 k_{j_l}$, for each $l \in \{ 1,2\}$.

By combining the three cases, we see that for all $i \in V(\Gamma)$, there
exists a line 
  bundle $L_i \stackrel{\psi_i}{\longrightarrow} E_i$ of degree $-e_i$ on
  $E_i$ and a meromorphic 2-form $\alpha_i$ of weight $r$ on it, such that:

   $\bullet$ $\div(\alpha_i)=  r \cdot k_i \: E_i  +  \sum_{j \in V_i} r \cdot k_j
       \: \psi_i^*(E_{i,j})$;  

   $\bullet$ for each point $A_{i,j,l}\in E_{i,j}$,  one has a system
       of local holomorphic coordinates $(q_{i,j,l}, f_{i,j,l})$ centered at
       $A_{i,j,l}$  on the total 
       space $L_i$, such that $q_{i,j,l}$ is constant on the fibers of
       $\psi_i$, that $f_{i,j,l}=0$ defines the $0$-section $E_i \subset L_i$ 
       and moreover: 
        \begin{equation} \label{formnorm2}
             \alpha_i= (f_{i,j,l})^{ r \cdot k_i} (q_{i,j,l})^{r \cdot k_j}  \:
             (df_{i,j,l} \wedge dq_{i,j,l})^r.  
         \end{equation}

Then, we plumb neighborhoods of the 0-sections of these bundles by identifying 
neighborhoods of the points $A_{i,j,l}$ and $A_{j,i,l}$, in such a way that the
forms $\alpha_i$ and $\alpha_j$ glue. Again, this is possible because we have
the simple normal forms (\ref{formnorm2}), and because no two
adjacent vertices $i$ and $j$ of $\Gamma$ satisfy $k_i = k_j = -1$, a fact
ascertained by the Corollary \ref{corgen} of Veys' Theorem \ref{core} 
(the analog of equation (\ref{coef}) is now the equation 
$\lambda_{i,j}^{r(k_j +1)} = (- 1)^r \lambda_{j,i}^{r(k_i  + 1)}$). 
After contracting the union of the 0-sections of the bundles $L_i$, one gets a
singularity $(X,0)$ such that the 2-form of weight $r$ obtained by gluing the
forms $\alpha_i$ gives
a holomorphic trivialization of $K_{X \setminus 0}^r$.
 The theorem is proved.    \hfill $\Box$
\medskip

The analog of Remark \ref{spesing} holds obviously also in this more general
context.  
\medskip 

To conclude, we mention that the problems of characterizing the topological
types  
of isolated complete intersections, hypersurfaces and smoothable isolated
surface  singularities are still open. It would also be interesting to determine whether 
the classes of $\Q$-Gorenstein singularities constructed in this paper are extremal among 
the $\Q$-Gorenstein singularities of the same topological type, for some analytical invariant.

\medskip
\emph{Acknowledgments}

I decided to do this research after collaborating with Jos\'e Seade 
on \cite{PS 09}. I am extremely grateful to him for his kind encouragements 
and for his remarks on a previous version of this paper. 
I am also grateful to Jonathan Wahl for having told me the problem and
for providing several remarks on earlier versions. Moreover, 
I am grateful both to him and to Bernard Teissier for having 
asked if I could extend Theorem  \ref{numeqgor} to $\Q$-Gorenstein
singularities,  and to Wim Veys for having told me about his paper \cite{V 04}.

\medskip

\bibliographystyle{plain}

Patrick Popescu-Pampu,

Universit\' e Lille 1, 
UFR de Math\' ematiques

B\^atiment M2, 
Cit\' e Scientifique

59655, Villeneuve d'Ascq Cedex

France.

\medskip
 patrick.popescu@math.univ-lille1.fr

\end{document}